\documentclass[12pt]{article}
\usepackage{amsfonts,amssymb,epsfig}
\usepackage{amsmath}

\date{}
\begin{document}
\newtheorem{df}{Definition}
\newtheorem{thm}{Theorem}
\newtheorem{lm}{Lemma}
\newtheorem{pr}{Proposition}
\newtheorem{co}{Corollary}
\newtheorem{re}{Remark}
\newtheorem{note}{Note}
\newtheorem{claim}{Claim}
\newtheorem{problem}{Problem}

\def\R{{\mathbb R}}

\def\E{\mathbb{E}}
\def\calF{{\cal F}}
\def\N{\mathbb{N}}
\def\calN{{\cal N}}
\def\calH{{\cal H}}
\def\n{\nu}
\def\a{\alpha}
\def\d{\delta}
\def\t{\theta}
\def\e{\varepsilon}
\def\t{\theta}
\def\pf{ \noindent {\bf Proof: \  }}
\def\trace{\rm trace}
\newcommand{\qed}{\hfill\vrule height6pt
width6pt depth0pt}
\def\endpf{\qed \medskip} \def\colon{{:}\;}
\setcounter{footnote}{0}

\def\Lip{{\rm Lip}}

\renewcommand{\qed}{\hfill\vrule height6pt  width6pt depth0pt}

\title{Tight embedding of subspaces of $L_p$ in $\ell_p^n$ for even $p$
\thanks {AMS subject classification: 46B07
}}

\author{Gideon Schechtman\thanks{Supported by the
Israel Science Foundation }}
\date{September 5, 2010}
\maketitle

\begin{abstract}

Using a recent result of Batson, Spielman and Srivastava, We obtain a tight estimate on the dimension of $\ell_p^n$, $p$ an even integer, needed to almost isometrically contain all $k$-dimensional subspaces of $L_p$.

\end{abstract}

In a recent paper \cite{bss} Batson, Spielman and Srivastava introduced a new method for sparsification of graphs which already proved to have functional analytical applications. Here we bring one more such application. Improving over a result of \cite{blm} (or see \cite{js} for a survey on this and related results), we show that for even $p$ and for $k$ of order $n^{2/p}$ any $k$ dimensional subspace of $L_p$ nicely embeds into $\ell_p^n$. This removes a $\log$ factor from the previously known estimate. The result in Theorem \ref{thm:embedding} is actually sharper than stated here and gives the best possible result in several respects, in particular in the dependence of $k$ on $n$.

The theorem of \cite{bss} we shall use is not specifically stated in \cite{bss}, but is stated as Theorem 1.6 in Srivastava's thesis \cite{sr}:

\begin{thm}{[BSS]}\label{thm:bss}
Suppose $0 < \e < 1$ and $A = \sum_{i=1}^mv_iv_i^T$
are given, with $v_i$ column vectors in $\R^k$. Then there are nonnegative weights $\{s_i\}_{i=1}^m$, at most $\lceil k/\e^2\rceil$
of which are nonzero, such that, putting $\tilde A=\sum_{i=1}^ms_iv_iv_i^T$,
\[
(1-\e)^{-2}x^TAx\le x^T\tilde Ax\le (1+\e)^{2}x^TAx
\]
for all $x\in \R^k$.
\end{thm}

\begin{co} \label{co:sampling} Let $X$ be a $k$-dimensional subspace of $\ell_2^m$ and let $0<\e<1$. Then there is a set $\sigma\subset\{1,2,\dots,m\}$ of cardinality $n\le C\e^{-2}k$ ($C$ an absolute constant) and positive weights $\{s_i\}_{i\in\sigma}$ such that
\[
\|x\|_2\le (\sum_{i\in\sigma}s_ix^2(i))^{1/2}\le(1+\e)\|x\|_2
\]
for all $x=(x(1),x(2),\dots,x(m))\in X$.
\end{co}

\pf Let $u_1,u_2,\dots,u_k$ be an orthonormal basis for $X$; $u_j=(u_{1,j},u_{2,j},\dots,u_{m,j})$, $j=1,\dots,k$.
Put $v_i^T=(u_{i,1},u_{i,2},\dots,u_{i,k})$, $i=1,\dots,m$. Then $\sum_{i=1}^m v_iv_i^T=I_k$, the $k\times k$  identity matrix. Let $s_i$ be the weights given by Theorem \ref{thm:bss} (and $\sigma$ their support). Then, for all $x=\sum_{i=1}^k a_iu_i\in X$,
\[
(1-\e)^{-2}\|x\|_2=a^T\sum_{i=1}^m v_iv_i^Ta\le a^T\sum_{i=1}^m s_iv_iv_i^Ta\le (1+\e)^{2}\|x\|_2.
\]
Finally, notice that, for each $i=1,\dots,m$, $a^T v_iv_i^Ta=x(i)^2$, the square of the $i$-th coordinate of $x$. Thus,
\[
a^T\sum_{i=1}^m s_iv_iv_i^Ta=\sum_{i=1}^ms_ix(i)^2.
\]
\endpf

We first prove a simpler version of the main result.

\begin{pr}\label{pr:embedding}
Let $X$ is a $k$ dimensional subspace of $L_p$ for some even $p$ and let $0<\e<1$. Then $X$ \ $(1+\e)$-embeds in $\ell_p^n$ for $n=O((\e p)^{-2}k^{p/2})$.
\end{pr}

\pf Assume as we may that $X$ is a $k$ dimensional subspace of $\ell_p^m$ for some finite $m$. Consider the set of all vectors which are coordinatewise products of $p/2$ vectors from $X$; i.e, of the form
\[
(x_1(1)x_2(1)\dots x_{p/2}(1), x_1(2)x_2(2)\dots x_{p/2}(2),\dots,x_1(m)x_2(m)\dots x_{p/2}(m))
\]
where $x_j=(x_j(1),x_j(2),\dots,x_j(m))$, $j=1,2,\dots,p/2$, are elements of $X$. We shall denote the vector above as $x_1\cdot x_2\cdot\dots\cdot x_{p/2}$. The span of this set in $\R^m$, which we denote by $X^{p/2}$, is clearly a linear space of dimension at most $k^{p/2}$. Consequently, by Corollary \ref{co:sampling}, there is a set $\sigma\subset \{1,2,\dots,m\}$ of cardinality at most $C(\e p)^{-2}k^{p/2}$ and weights $\{s_i\}_{i\in\sigma}$ such that
\[
\|y\|_2\le (\sum_{i\in\sigma}s_iy^2(i))^{1/2}\le(1+\frac{\e p}{4})\|y\|_2
 \]
for all $y\in X^{p/2}$. Restricting to $y$-s of the form \[
y=(x(1)^{p/2},x(2)^{p/2},\dots,x(m)^{p/2})
\]
with $x=(x(1),x(2),\dots, x(m))\in X$, we get
\[
\|x\|_p^{p/2}\le (\sum_{i\in\sigma}s_ix^p(i))^{1/2}\le(1+\frac{\e p}{4})\|x\|_p^{p/2}.
\]
Raising these inequalities to the power $2/p$ gives the result.
\endpf

We now state and prove the main result.

\begin{thm}\label{thm:embedding}
Let $X$ be a $k$ dimensional subspace of $L_p$ for some even $p\le k$ and let $0<\e<1$. Then $X$ \ $(1+\e)$-embeds in $\ell_p^n$ for $n=O(\e^{-2} (\frac{10 k}{p})^{p/2})$. Equivalently, for some universal $c>0$, for any $n$ and any $k\le c\e^{4/p}pn^{2/p}$, any $k$-dimensional subspace of $L_p$  \ $(1+\e)$-embeds in $\ell_p^n$.
\end{thm}

\pf The only change from the previous proof is a better estimate of the dimension of the auxiliary subspace involved. An examination of the proof above shows that it is enough to apply Corollary \ref{co:sampling} to any subspace containing all the vectors $x^{p/2}=x\cdot x\cdot\dots\cdot x$ ($p/2$ times), $x=(x(1),\dots,x(m))\in X$. The smallest such subspace is the space of degree $p/2$ homogeneous polynomials in elements of $X$. Its basis is the set of monomials of the form $u_{j_1}^{p_1}\cdot u_{j_2}^{p_2}\cdot\dots\cdot u_{j_l}^{p_l}$ with $p_1+\dots + p_l=p/2$, where $u_1,\dots,u_k$ is a basis for $X$.   The dimension of this space, which is the number of such monomials, is ${k+p/2-1\choose p/2} \le (\frac{10 k}{p})^{p/2}$.
\endpf

Next we remark on the estimate  $k\le c\e^{4/p}pn^{2/p}$.\hfill\break
 This estimate improves (unfortunately, only for even $p$) over the known estimates (the best of which is in \cite{blm}) by removing a $\log n$ factor that was presented in the best estimate till now. Also, the $p$ in front of the $n^{2/p}$ is a nice feature. It is known that the dependence of $k$ on $p$ and $n$ in this estimate is best possible even if one restrict to subspaces of $L_p$ isometric to $\ell_2^k$ (see \cite{bdgjn}). Actually the result above indicates that $\ell_2^k$ are the ``worse" subspaces.
 
 As for the dependence on $\e$, the published proofs gives at best quadratic dependence while here we get a quadratic dependence for $p=4$ and better ones as $p$ grows. For the special case of $X=\ell_2^k$ and $p=4$ a better result is known: One can embed $\ell_2^k$ isometrically into $\ell_4^{4n^2}$ (\cite{ko}). But for $p=6,8,\dots$ we get better result here even for this special case than what was previously known (The best I knew was a linear dependence on $\e$ - this is unpublished. Here we get $\e^{2/3}$ for $p=6$ and better for larger $p$-s.) It is not clear that this is the best possible dependence on $\e$, but note that for some combinations of $\e$ and $p$ and in particular for every $\e$ and large enough $p$ ($p>c\log 1/\e$) the dependence on $\e$ becomes a constant and, up to universal constants, we get the best possible result with respect to all parameters.

\vfill\eject

\noindent Gideon Schechtman\newline Department of
Mathematics\newline Weizmann Institute of Science\newline Rehovot,
Israel\newline E-mail: gideon.schechtman@weizmann.ac.il

\end{document}